\documentclass[reqno,11pt]{amsart} \numberwithin{equation}{section}
\input{amssym.def}
\input{amssym.tex}
\usepackage{graphicx}
\begin{document}
\title[Convergence theorems by extragradient method]
{Convergence theorems by extragradient method in Banach spaces}
\author[Z. Jouymandi and F. Moradlou]{Zeynab Jouymandi$^1$ and Fridoun Moradlou$^2$}
\address{\indent $^{1,2}$ Department of Mathematics
\newline \indent Sahand University of Technology
\newline \indent Tabriz, Iran}
\email{\rm $^1$ z\_jouymandi@sut.ac.ir \& z.juymandi@gmail.com}
\email{\rm $^2$ moradlou@sut.ac.ir \& fridoun.moradlou@gmail.com}
\thanks{}
\begin{abstract}
In this paper, using generalized metric projection, we propose a new extragradient method for finding a
common element of the solutions set of a generalized equilibrium
problem and a variational inequality for an
$\alpha$-inverse-strongly monotone operator and fixed points of two
relatively nonexpansive mappings in Banach spaces. We prove strong
convergence theorems by this method under suitable conditions. A
numerical example is given to illustrate the usability of our
results.
\end{abstract}
\keywords{Generalized equilibrium problem, Generalized metric
projection, Relatively nonexpansive mapping, Variational
inequality, Strong convergence} \subjclass[2010]{Primary 65K10,
47H09,  47J05, 47J25.} \maketitle

\newtheorem{df}{Definition}[section]
    \newtheorem{rk}[df]{Remark}
 \newtheorem{lm}[df]{Lemma}
  \newtheorem{thm}[df]{Theorem}
  \newtheorem{exa}[df]{Example}
  \newtheorem{cor}[df]{Corollary}

\section{Introduction}
Let $E$ be a real Banach space and $E^{*}$ be the dual of $E$. Let
$C$ be a closed convex subset of $E$. In this paper, we concerned
with the following Variational Inequality ($VI$), which consists
in finding a point $u\in C$ such that
\begin{eqnarray*}\label{vthoe1}
\langle Au,y-u\rangle\geq 0, \qquad \forall \ y\in C,
\end{eqnarray*}
where $A:C\rightarrow E^{*}$ is a given mapping and
$\langle.,.\rangle$ denotes the generalized duality pairing. The
solution set of ($VI$) denoted by $SOL(C,A)$.
\par
 Let
$A:C\rightarrow E^{*}$ be a nonlinear mapping and $f:C\times
C\rightarrow \mathbb{R}$ be a bifunction. We consider the
following generalized equilibrium problem of finding $u\in C$ such
that
\begin{equation}\label{vthoe2}
f(u,y)+\langle Au,y-u\rangle\geq 0, \qquad \forall \ y\in C.
\end{equation}
The solutions set of (\ref{vthoe2}) is denoted by $GEP(f,A)$,
i.e.,
\begin{eqnarray*}
    GEP(f,A)=\{u\in C: f(u,y)+\langle Au,y-u\rangle\geq 0,  \ \forall \ y\in C\}.
\end{eqnarray*}
In the case of $A\equiv 0$, problem (\ref{vthoe2}) is equivalent
to finding $u\in C$
 such that $f(u,y)\geq 0$, for all $y\in C$, which is called the equilibrium problem.
 The set of its solutions is denoted by $EP(f)$. In the case of $f\equiv0$,
 the problem (\ref{vthoe2}) reduces to $(VI)$.
 \par
 A mapping $T:C\rightarrow C$ is called nonexpansive if
$$\Vert T(x)-T(y) \Vert\leqslant \Vert x-y\Vert,$$ for all $x, y\in
C$. The set of fixed points of $T$ is the set $F(T):=\{x\in C:
Tx=x\}$. An operator $A:C\rightarrow E^{*}$ is called monotone if
\[\langle Ax-Ay,x-y\rangle \geq 0,\]
for all $x, y\in C$. Also, it is called $\alpha$-inverse-strongly
monotone if there exists a constant $\alpha>0$ such that \[
\langle Ax-Ay,x-y \rangle \geqslant \alpha \Vert Ax-Ay
\Vert^{2},\]for all $x,y\in C$. A monotone operator $A$ is said to
be maximal if its graph $G(A)=\{(x,Ax): x\in D(A)\}$ is not
contained in the graph of any other monotone operator. It is clear
that a monotone operator $A$ is maximal if and only if, for any
$(x,x^{*})\in E\times E^{*}$, if $\langle
x-y,x^{*}-y^{*}\rangle\geq0$ for all $(y,y^{*})\in G(A)$, then it
follows that $x^{*}=Ax$.
\par
In $1976$, Korpelevich \cite{korp} proposed a new algorithm for
solving the ($VI$) in the Euclidean space which is known as
extragradient method. Putting $x^{0}\in H$ arbitrarily, she present
her algorithm as follows:
$$\begin{cases}
 y^{k}:=P_{C}(x^{k}-\tau Ax^{k})\\
x^{k+1}:=P_{C}(x^{k}-\tau Ay^{k}),
\end{cases}$$
where $\tau$ is a positive number and $P_{C}$ denotes Euclidean
least distance projection of $H$ onto $C$.
\par
In $2008$, Plubtieng and Punpaeng \cite{Plub} have introduced the
following iteration process for finding a common element of
solutions set of a ($VI$) for an $\alpha$-inverse-strongly
monotone operator $A$, the set of solutions of an equilibrium
problem and the set of fixed points of a nonexpansive mapping $S$
with $\Omega=SOL(C,A)\cap EP(f)\cap F(S)\neq\emptyset$, in a
Hilbert space. Let $\{x^{k}\}$ be a sequence generated by
$$\begin{cases}
x^{1}\in C,\\
u^{k}\in C\ \ \ \text{such that}\ \ \ f(u^{k},y)+\frac{1}{r^{k}}\langle y-u^{k},u^{k}-x^{k}\rangle\geq0,\quad\forall\: y\in C,\\
y^{k}:=P_{C}(u^{k}-\tau Au^{k}),\\
x^{k+1}:=\alpha^{k}x^{1}+\beta^{k}x^{k}+\gamma^{k}SP_{C}(y^{k}-\lambda^{k}Ay^{k}),\quad\forall\:
k\geq 1,
\end{cases}$$
where $P_{C}$ denotes metric projection of $H$ onto $C$,
$\{\alpha^{k}\}$, $\{\beta^{k}\}$ and $\{\gamma^{k}\}$ are
sequences in $[0,1)$ and $\{\lambda^{k}\}$ is a sequence in
$[0,2\alpha]$. Under suitable conditions, they have proved
$\{x^{k}\}$ converges strongly to $P_{\Omega}x^{1}$.
\par
 Very recently, Qin et al.
\cite{Qin} have introduced the following iteration process for two
relatively nonexpansive mappings. The sequence $\{x^{k}\}$
generated by
$$\begin{cases}
u^{1}\in C,\\
x^{k}\in C\ \ \ \text{such that}\ \ \ f(x^{k},y)+\frac{1}{r^{k}}\langle y-x^{k},Jx^{k}-Ju^{k}\rangle\geq0,\quad\forall\: y\in C,\\
u^{k
+1}:=J^{-1}(\alpha^{k}Jx^{k}+\beta^{k}JTx^{k}+\gamma^{k}JSx^{k}),\quad\forall\:
k\geq 1,
\end{cases}$$
converges weakly to $\nu\in \Omega=F(T)\cap F(S)\cap EP(f)$, where
$\{\alpha^{k}\}$, $\{\beta^{k}\}$ and $\{\gamma^{k}\}$ satisfy
suitable conditions,
$\nu=\lim\limits_{k\rightarrow\infty}\Pi_{\Omega}x^{k}$ and
$\Pi_{\Omega}$ denotes generalized metric projection in  Banach
spaces which is an analogue of the metric projection in Hilbert
spaces.
\par
The extragradient algorithm is well known because of its
efficiency in numerical tests. Therefore, in recent years, many
authors have used extragradient method for finding a common
element of solutions set of a ($VI$), the set of solutions of an
equilibrium problem and the set of fixed points of a nonexpansive
or a relatively nonexpansive mapping in the framework of Hilbert
spaces and Banach spaces, see for instance \cite{tak,Ta} and the
references there in. In all of these methods, authors have proved weak convergence of generated sequences.
 \par
In this paper, motivated Plubtieng and Punpaeng \cite{Plub} and Qin
et al.\cite{Qin}, we propose a new extragradient algorithm by using
generalized metric projection. Using this algorithm, we
prove strong convergence theorems under suitable conditions.
\section{Preliminaries}
 Let $E$ be a real Banach space and $E^{*}$ be the dual of $E$.
 The normalized duality mapping $J$ from $E$ to $2^{E^{*}}$ is defined by
\[Jx=\lbrace x^{*}\in E^{*} :\langle x , x^{*} \rangle=\Vert x\Vert^{2}=\Vert x^{*}\Vert^{2}\rbrace, \qquad \forall  \ x\in E.\]
Also, the strong convergence and the weak convergence of a
sequence $\{x^{k}\}$ to $x$ in $E$ are denoted by
$x^{k}\rightarrow x$ and $x^{k}\rightharpoonup x$, respectively.
\par
Let $S(E)$ be the unite sphere centered at the origin of $E$. A
Banach space $E$ is strictly convex if $\|\frac{x+y}{2}\|<1$,
whenever $x, y \in S(E)$ and $x\neq y$. Modulus of convexity of
$E$ is defined by $$\delta_{E}(\epsilon)= \inf
\{1-\frac{1}{2}\|(x+y)\|:\ \|x\|, \|y\|\leq 1, \ \|x-y\|\geq
\epsilon\},$$ for all $\epsilon\in [0, 2]$. Also, $E$ is said to be
uniformly convex if $\delta_{E}(0)=0$ and
$\delta_{E}(\epsilon)>0,$ for all $0<\epsilon\leq2$. Let $p$ be a
fixed real number with $p\geq2$. A Banach space $E$ is said to be
$p$-uniformly convex \cite{kato} if there exists a constant $c>0$
such that $\delta_{E}\geq c\epsilon^{p}$ for all $\epsilon\in [0,
2]$. The Banach space $E$ is called smooth if the limit
 \begin{equation}\label{eq1}
 \lim_{t\rightarrow0}\frac{\|x+ty\|-\|x\|}{t},
 \end{equation}
 exists for all $x,y\in S(E)$. It is also said to be uniformly smooth if the limit
  $(\ref{eq1})$ is attained uniformly for all $x,y\in S(E)$. Every uniformly smooth Banach space $E$ is smooth.
If a Banach space $E$ uniformly convex, then $E$ is reflexive and
strictly convex. For more details see \cite{Agarwal,Takahashi}.
  \par
  Some properties of the normalized duality mapping $J$ are listed in the following:
  \begin{enumerate}
\item For every $x\in E$, $Jx$ is nonempty closed convex and bounded subset of $E^*$.
\item  If $E$ is smooth or $E^*$ is strictly convex, then $J$ is single-valued.
\item If $E$ is strictly convex, then $J$ is one-one, i.e., if $x\neq y$ then $ Jx\cap Jy=\phi$.
\item If $E$ is reflexive, then $J$ is onto.
\item If $E$ is strictly convex, then $J$ is strictly monotone, that is, $$\langle x-y, Jx-Jy\rangle>0,$$
for all  $x, y\in E$ such that $x\neq y$.
\item If $E$ is smooth and reflexive, then $J$ is norm-to-weak$^{*}$ continuous, that is, $Jx^{k}\rightharpoonup^{*}Jx$ whenever $x^{k}\rightarrow x$.
\item If $E$ is smooth, strictly convex  and reflexive  and $J^*:E^*\rightarrow 2^E$ is the normalized duality mapping on $E^*$, then $J^{-1}=J^*$, $JJ^*=I_{E^*}$ and $J^*J=I_{E}$, where $I_{E}$ and $I_{E^*}$ are the identity mapping on $E$ and $E^*$, respectively.
\item If $E$ is  uniformly convex and uniformly smooth, then $J$ is uniformly norm-to-norm continuous on bounded sets of $E$ and $J^{-1}=J^*$ is also uniformly norm-to-norm continuous on bounded sets of $E^*$, i.e., for $\varepsilon>0$ and $M>0$, there is a $\delta>0$ such that
\begin{equation}\label{eq30}
\|x\|\leq M, \ \|y\|\leq M \ \ \text{and} \ \ \|x-y\|<\delta \ \
\Rightarrow \ \ \|Jx-Jy\|<\varepsilon,
\end{equation}
\begin{equation}\label{eq18}
\|x^{*}\|\leq M, \ \|y^{*}\|\leq M \ \ \text{and} \ \
\|x^{*}-y^{*}\|<\delta \ \ \Rightarrow \ \
\|J^{-1}x^{*}-J^{-1}y^{*}\|<\varepsilon.
\end{equation}
\end{enumerate}
\par
Let $E$ be a smooth Banach space, the function $\phi: E\times
E\rightarrow \mathbb{R}$ is defined by \[\phi(x, y)=\Vert
x\Vert^{2}-2\langle x,Jy \rangle+\Vert y\Vert^{2},\] for all $x,
y\in E$. It is clear from definition of $\phi$ that for all $x, y,
z, w\in E$,
\begin{equation}\label{eq14}
(\Vert x\Vert -\Vert y\Vert)^{2} \leq \phi(x,y)\leq (\Vert
x\Vert+\Vert y\Vert)^{2},
\end{equation}
\begin{equation}\label{eq36}
\phi(x,y)=\phi(x,z)+\phi(z,y)+2\langle x-z,Jz-Jy\rangle,
\end{equation}
\begin{equation}\label{eq500}
\phi(x,y)=\langle x,Jx-Jy\rangle+\langle y-x,Jy\rangle\leq\|x\|\|Jx-Jy\|+\|x-y\|\|y\|.
\end{equation}
Also, the function $V:E\times E^{*}\rightarrow \mathbb{R}$ is
defined by $V(x, x^{*})=\| x\|^{2}-2<x, x^{*}>+\| x^{*} \|^{2}$,
for all $x\in E$ and all $x^{*}\in E$. That is,
$V(x,x^{*})=\phi(x,J^{-1}x^{*})$ for all $x\in E$ and all $x\in
E^{*}$. It is well known that, if $E$ is a reflexive strictly
convex and smooth Banach space with $E^{*}$ as its dual, then
\begin{equation}\label{eq5}
 V(x, x^{*})+2\langle J^{-1}x^{*}-x,y^{*}\rangle \leq V(x, x^{*}+y^{*}),
 \end{equation}
 for all $x\in E$ and all $x^{*},y^{*}\in E^{*}$ \cite{Rockfellar}.
\par
An operator $A:C\rightarrow E^{*}$ is called hemicontinuous, if for all $x, y\in C$, the mapping $T$ of $[0,1]$ into $E^{*}$ defined by $T(t)=A(tx+(1-t)y)$ is continuous with respect to the weak$^{*}$ topology of $E^{*}$.
\par
Let $C$ be a closed convex subset of a smooth Banach space $E$ and
$T:C\rightarrow C$ be a mapping. A point $p$ in $C$ is said to be
an asymptotic fixed point of $T$ if $C$ contains a sequence
$\{x^{k}\}$ which converges weakly to $p$ such that
$\lim\limits_{k\rightarrow\infty}(Tx^{k}-x^{k})=0$. The set of
asymptotic fixed points of $T$ will be denoted by $\hat{F}(T)$. A
mapping $T:C\rightarrow C$ is called relatively nonexpansive if
$\hat{F}(T)=F(T)$ and $\phi(p,Tx)\leq\phi(p,x)$ for all $x\in C$
and all $p\in F(T)$. The asymptotic behavior of relatively
nonexpansive mappings was studied in \cite{But}.
The mapping $T$ is said to be relatively quasi-nonexpansive if $F(T)\neq\emptyset$ and $\phi(p,Tx)\leq\phi(p,x)$ for all $x\in C$ and all $p\in F(T)$. The class of relatively quasi-nonexpansive mapping is broader than the class of relatively nonexpansive mappings which requires $\hat{F}(T)=F(T)$.\\
It is well known that, if $E$ is a strictly convex and smooth
Banach space, $C$ is a nonempty closed convex subset of $E$ and
$T:C\rightarrow C$ is a relatively quasi-nonexpansive mapping,
then $F(T)$ is a closed convex subset of $C$ \cite{Q}.
\begin{lm}\cite{Alber}\label{eq4}
Let $C$ be a nonempty closed subset of a smooth, reflexive and strictly
convex Banach space $E$ such that $\Pi_{C}$ from $E$ onto $C$ be a generalized metric projection and let $(z,x)\in
C\times E$. Then the following hold:
\begin{enumerate}
\item $z=\Pi_{C}x $ if and only if $\langle y-z,Jx-Jz\rangle\leq 0$ for all $y\in C$,
\item $\phi(z,\Pi_{C}x)+\phi(\Pi_{C}x,x)\leq\phi(z,x),$
\item $\phi(\Pi_{C}x,x)=\min_{y\in C}\phi(y,x)$
\end{enumerate}
\end{lm}
We need the following lemmas for the proof of our main results in
next section.
\begin{lm}\cite{Agarwal}
Let $E$ be a topological space and $f:E\rightarrow
(-\infty,\infty]$ be a function. Then the following statements are
equivalent:
\begin{enumerate}
\item $f$ is lower semicontinuous.
\item For each $\alpha\in \mathbb{R}$, the level set $\{x\in E: \ \ f(x)\leq\alpha\}$ is closed.
\item The epigraph of the function $f$, $\{(x,\alpha)\in E\times \mathbb{R}: \ \ f(x)\leq\alpha\}$ is closed.
\end{enumerate}
\end{lm}
\begin{lm}\cite{Agarwal}\label{eq23}
Let $C$ be nonempty closed convex subset of a Banach space $E$ and
$f:E\rightarrow (-\infty,\infty]$ be a convex function. Then $f$
is lower semicontinuous in the norm topology if and only if $f$ is
lower semicontinuous in the weak topology.
\end{lm}
\begin{lm}\cite{Yao}\label{eq6}
Let $E$ be a $2$-uniformly convex and smooth Banach space. Then,
for all $x, y\in E$, we have
$$\|x-y\|\leq\frac{2}{c^{2}}\|Jx-Jy\|,$$ where $\frac{1}{c}(0\leq
c\leq 1)$ is the $2$-uniformly convex constant of $E$.
\end{lm}
\begin{lm}\label{kamimra}\cite{Kamimura}
Let $E$ be a uniformly convex Banach space and let $r>0$. Then
there exists a strictly increasing, continuous and convex function
$g:[0, 2r]\rightarrow [0, \infty)$, $g(0)=0$ such that
$$g(\|x-y\|)\leq\phi(x,y),$$ for all $x,y\in B_{r}(0)=\{z\in
E:\|z\|\leq r\}$.
\end{lm}
\begin{lm}\cite{Cho}\label{eq15}
Let $E$ be a uniformly convex Banach space. Then there exists a
continuous strictly increasing convex function $g:[0,
2r]\rightarrow [0, \infty)$, $g(0)=0$ such that $$\|\lambda x+\mu
y+\gamma z\|^{2}\leq
\lambda\|x\|^{2}+\mu\|y\|^{2}+\gamma\|z\|^{2}-\lambda\mu
g(\|x-y\|),$$ for all $x,y,z\in B_{r}(0)=\{z\in E:\|z\|\leq r\}$
and all $\lambda, \mu, \gamma\in [0, 1]$ with
$\lambda+\mu+\gamma=1$.
\end{lm}
\begin{lm}\cite{Kamimura}\label{eq22}
Let $E$ be a uniformly convex and smooth Banach space and let
$\{x^{k}\}$ and $\{y^{k}\}$ be two sequences in $E$. If
$\phi(x^{k}, y^{k})\rightarrow 0$ and either $\{x^{k}\}$ or
$\{y^{k}\}$ is bounded, then $x^{k}-y^{k}\rightarrow 0.$
\end{lm}
We denote by $N_{C}(\nu)$ the normal cone for $C$ at a point
$\nu\in C$, that is \[N_{C}(\nu):=\{x^{*}\in E^{*}: \langle
\nu-y,x^{*}\rangle\geq0,\ \forall  \ y\in C\}.\]
\begin{lm}\cite{Rockfellar}\label{eq34}
Let $C$ be a nonempty closed convex subset  of a Banach space $E$
and let $A$ be a monotone and hemicontinuous operator of $C$ into
$E^{*}$ with $C=D(A)$. Let $B\subset E\times E^{*}$ be an operator
define as follows:
$$B\nu=\begin{cases}A\nu+N_{C}(\nu), \ \ \nu\in C,\\
\emptyset, \ \ \ \ \ \ \ \ \ \ \ \ \  \ \ \nu\not \in C.
\end{cases}$$
Then $B$ is maximal monotone and $B^{-1}(0)= SOL(A,C)$.
\end{lm}
\par
For solving the generalized equilibrium problem, we assume that
$f:C\times C\longrightarrow\mathbb{R}$ satisfies the following
conditions:
\begin{enumerate}
\item[(A1)] $f(x,x)=0$ for all $x\in C,$
\item[(A2)]$f$ is monotone, i.e., $f(x,y)+f(y,x)\leq0$ for all $x,y\in C,$
\item[(A3)]for each $x,y,z\in C,\ \
\lim\limits_{t\downarrow 0}f(t z+(1-t)x,y)\leq f(x,y),$
\item[(A4)]for each $x\in C,\;y\mapsto f(x,y)$ is convex and lower semicontinuous.
\end{enumerate}
\begin{lm}\cite{Liou}\label{eq13}
 Let $E$ be a smooth, strictly convex and reflexive Banach space and $C$ be a nonempty closed convex subset of $E$.
 Let $A:C\rightarrow E^{*}$ be an $\alpha$-inverse-strongly monotone operator,
 $f$ be a bifunction from $C\times C$ to $\mathbb{R}$ satisfying $(A1)-(A4)$ and
let $r >0$. Then for all $x\in E$, there exists $u\in C$ such that
$$f(u,y)+\langle Au,y-u\rangle+\frac{1}{r}\langle y-u,Ju-Jx\rangle\geq 0,\; \ \forall  \ y\in C,$$
if $E$ is additionally uniformly smooth and $K_{r}:E\rightarrow C$
is defined as
\begin{equation}\label{eq70}
K_{r}x = \{u\in C : f(u,y)+\langle
Au,y-u\rangle+\frac{1}{r}\langle y-u,Ju-Jx\rangle\geq0,\; \
\forall  \ y\in C\}.
\end{equation}
Then, the following statements hold:
    \begin{enumerate}
      \item[(i)]$K_r$ is singel-valued,
      \item[(ii)]$K_r$ is firmly nonexpansive, i.e., for all $x,y\in E$,
      $$\langle K_{r}x-K_{r}y, JK_{r}x-JK_{r}y\rangle\leq\langle K_{r}x-K_{r}y,Jx-Jy\rangle,$$
      \item[(iii)]$F(K_{r})=\hat{F}(K_{r})=GEP(f,A)$,
      \item[(iv)]$GEP(f,A)$ is closed and convex,
      \item[ (v) ] $\phi(p,K_{r}x)+\phi(K_{r}x,x)\leq \phi(p,x), \ \forall \ p\in F(K_{r}).$
    \end{enumerate}
    \end{lm}
\section{Main results}
Now, we present a new extragradient algorithm for finding a solution of the ($VI$)
which is also the common element of the set of solutions of a
generalized equilibrium problem and the set of fixed points of two
relatively nonexpansive mappings.
\begin{thm}\label{eq25}
Let $C$ be a nonempty closed convex subset of a $2$-uniformly
convex, uniformly smooth Banach space $E$. Assume that $f:C\times
C\longrightarrow \mathbb{R}$ is a bifunction which satisfies
conditions $(A1)-(A4)$. Let $A:C\rightarrow E^{*}$ be an
$\alpha$-inverse strongly monotone operator and $T,S:C\rightarrow
C$ be two relatively nonexpansive mappings such that
\[\Omega:=SOL(C,A)\cap GEP(f,A)\cap F(T)\cap F(S)\neq\emptyset,\]
and $\Vert Ax\Vert\leq\Vert Ax-Au\Vert   \  \text{for all}\ x\in
C\ \text{and all} \ u\in  \Omega$. Assume that $\Pi_{C}$ is the
generalized metric projection from $E$ onto $C$. Let
$\{x^{k}\}$ be a sequence generated by $x^{1}\in C$ and
\begin{equation}\label{eq3}
\begin{cases}
 u^{k}\in C\;\,\,\,\text{s.t.}\,\,\, f(u^{k},y)+\langle Au^{k},y-u^{k}\rangle +\frac{1}{r^{k}}\langle y-u^{k},Ju^{k}-Jx^{k}\rangle\geq0,\quad\forall\: y\in C,&\\
  y^{k}:=\Pi_{C}J^{-1}(Jx^{k}-\tau Ax^{k}),&\\
 z^{k}:=\Pi_{C}J^{-1}(Ju^{k}-\tau Au^{k}),&\\
 x^{k+1}:=\Pi_{C}J^{-1}(\alpha^{k}Jx^{k}+\beta^{k}JTz^{k}+\gamma^{k}JSy^{k}).&
  \end{cases}
  \end{equation}
Furthermore, suppose that $\{\alpha^{k}\}$, $\{\beta^{k}\}$ and
$\{\gamma^{k}\}$ are three sequences in $[0,1]$ satisfying the
following conditions:
 \begin{enumerate}
      \item[(i)]$\alpha^{k}+\beta^{k}+\gamma^{k}=1$,
      \item[(ii)]$\liminf\limits_{k\rightarrow \infty}\alpha^{k}\beta^{k}>0 \quad \& \quad \liminf\limits_{k\rightarrow \infty}\alpha^{k}\gamma^{k}>0$,
      \item[(iii)]$\{r^{k}\}\subset [a,\infty)$ for some $a>0$,
      \item[(iv)]$0<\tau<\frac{c^{2}\alpha}{2}$, where $\frac{1}{c}\,\,(0< c\leq 1)$ is the $2$-uniformly convexity constant of $E$.
    \end{enumerate}
Then the sequences $\{x^{k}\}_{k=1}^{\infty}$,
$\{y^{k}\}_{k=1}^{\infty}$ and $\{z^{k}\}_{k=1}^{\infty}$
generated by (\ref{eq3}) converge strongly  to the some solution
$u^{*}\in \Omega$, where $u^{*}=\lim\limits_{k\rightarrow
\infty}\Pi_{\Omega}({x^{k}}).$
\end{thm}
{\it Proof} Let $u\in\Omega$, from Lemma \ref{eq4}, the definition
of function $V$ and inequality (\ref{eq5}), we get
\begin{equation}
\begin{aligned}
\phi(u,y^{k})&=\phi(u,\Pi_{C}J^{-1}(Jx^{k}-\tau Ax^{k}))\\
&\leq\phi(u,J^{-1}(Jx^{k}-\tau Ax^{k}))\\
&=V(u,(Jx^{k}-\tau Ax^{k}))\\
&\leq V(u,(Jx^{k}-\tau Ax^{k})+\tau Ax^{k})-2\langle J^{-1}(Jx^{k}-\tau Ax^{k})-u,\tau Ax^{k}\rangle \\
&=V(u,Jx^{k})-2\langle J^{-1}(Jx^{k}-\tau Ax^{k})-u,\tau Ax^{k}\rangle\\
&=\phi(u,x^{k})-2\tau\langle x^{k}-u,Ax^{k}\rangle+2\langle
J^{-1}(Jx^{k}-\tau Ax^{k})-x^{k},-\tau Ax^{k}\rangle.
\end{aligned}
\end{equation}
Since $A$ is an $\alpha$-inverse strongly monotone operator and
$u\in SOL(C,A)$, we have
\begin{equation}\label{eq7}
\begin{aligned}
-2\tau\langle x^{k}-u,Ax^{k}\rangle&=-2\tau\langle x^{k}-u,Ax^{k}-Au\rangle-2\tau\langle x^{k}-u,Au\rangle\\
&\leq -2\alpha\tau\|Ax^{k}-Au\|^{2}.
\end{aligned}
\end{equation}
 From Lemma \ref{eq6} and $\Vert Ax\Vert\leq\Vert
Ax-Au\Vert   \  \text{for all}\ x\in C\ \text{and all}\ u\in
\Omega$, we obtain
\begin{equation}\label{eq8}
\begin{aligned}
2\langle J^{-1}(Jx^{k}-\tau Ax^{k})&-x^{k},-\tau Ax^{k}\rangle\\
&=2\langle J^{-1}(Jx^{k}-\tau Ax^{k})-J^{-1}J(x^{k}),-\tau Ax^{k}\rangle\\
&\leq 2\|J(J^{-1}(Jx^{k}-\tau Ax^{k}))-J(J^{-1}Jx^{k})\|\|\tau Ax^{k}\|\\
&\leq \frac{4}{c^{2}}\tau^{2}\|Ax^{k}\|^{2}\\
&\leq \frac{4}{c^{2}}\tau^{2}\|Ax^{k}-Au\|^{2}.
\end{aligned}
\end{equation}
It follows from inequalities (\ref{eq7}), (\ref{eq8}) and
condition (iv) that
\begin{equation}\label{eq10}
\phi(u,y^{k})\leq\phi(u,x^{k})+2\tau(\frac{2\tau}{c^{2}}-\alpha)\|Ax^{k}-Au\|^{2}\leq\phi(u,x^{k}).
\end{equation}
In a similar way, we can conclude
\begin{equation}\label{eq21}
\phi(u,z^{k})\leq\phi(u,u^{k})+2\tau(\frac{2\tau}{c^{2}}-\alpha)\|Au^{k}-Au\|^{2}\leq\phi(u,u^{k}).
\end{equation}
From (\ref{eq3}) and condition (v) of Lemma \ref{eq13}, we have
\begin{equation}\label{eq32}
\phi(u,u^{k})=\phi(u,K_{r^{k}}x^{k})\leq\phi(u,x^{k}),
\end{equation}
hence, we conclude that
\begin{equation}\label{eq12}
\phi(u,z^{k})\leq\phi(u,x^{k}).
\end{equation}
By the convexity of $\|.\|^{2}$, the definition of $T,\ S$ and
inequalities (\ref{eq10}) and (\ref{eq12}), we obtain
\begin{equation}\label{eq20}
\begin{aligned}
\phi(u,x^{k+1})&\leq\phi(u,J^{-1}(\alpha^{k}Jx^{k}+\beta^{k}JTz^{k}+\gamma^{k}JSy^{k}))\\
&=\|u\|^{2}-2\alpha^{k}\langle u,Jx^{k}\rangle-2\beta^{k}\langle u,JTz^{k}\rangle-2\gamma^{k}\langle u,JSy^{k}\rangle\\
&\quad+\|\alpha^{k}Jx^{k}+\beta^{k}JTz^{k}+\gamma^{k}JSy^{k}\|^{2}\\
&\leq\|u\|^{2}-2\alpha^{k}\langle u,Jx^{k}\rangle-2\beta^{k}\langle u,JTz^{k}\rangle-2\gamma^{k}\langle u,JSy^{k}\rangle\\
&\quad+\alpha^{k}\|Jx^{k}\|^{2}+\beta^{k}\|JTz^{k}\|^{2}+\gamma^{k}\|JSy^{k}\|^{2}\\
&=\alpha^{k}\phi(u,x^{k})+\beta^{k}\phi(u,Tz^{k})+\gamma^{k}\phi(u,Sy^{k})\\
&\leq\alpha^{k}\phi(u,x^{k})+\beta^{k}\phi(u,z^{k})+\gamma^{k}\phi(u,y^{k})\\
&\leq\phi(u,x^{k}).
\end{aligned}
\end{equation}
This implies that $\lim\limits_{k\rightarrow\infty}\phi(u,x^{k})$
exists. Therefore $\{\phi(u,x^{k})\}$ is bounded. From inequality
(\ref{eq14}), we know that $\{x^{k}\}$ is bounded. Therefore, it
follows from inequalities (\ref{eq10}), (\ref{eq32}) and
(\ref{eq12}) that $\{y^{k}\}$, $\{u^{k}\}$ and $\{z^{k}\}$ are
also bounded. Let $r_{1}=\sup_{k\geq1}\{\|x^{k}\|,\|Tz^{k}\|\}$
and $r_{2}=\sup_{k\geq1}\{\|x^{k}\|,\|Sy^{k}\|\}$. So, by Lemma
\ref{eq15}, there exists a continuous, strictly increasing and
convex function $g_{1}:[0,2r_{1}]\rightarrow\mathbb{R}$ with
$g_{1}(0)=0$ such that for $u\in \Omega$, we get
\begin{equation*}
\begin{aligned}
\phi(u,x^{k+1})&\leq\|u\|^{2}-2\alpha^{k}\langle u,Jx^{k}\rangle-2\beta^{k}\langle u,JTz^{k}\rangle\\
&\quad-2\gamma^{k}\langle u,JSy^{k}\rangle+\|\alpha^{k}Jx^{k}+\beta^{k}JTz^{k}+\gamma^{k}JSy^{k}\|^{2}\\
&\leq\|u\|^{2}-2\alpha^{k}\langle u,Jx^{k}\rangle-2\beta^{k}\langle u,JTz^{k}\rangle-2\gamma^{k}\langle u,JSy^{k}\rangle\\
&\quad+\alpha^{k}\|Jx^{k}\|^{2}+\beta^{k}\|JTz^{k}\|^{2}+\gamma^{k}\|JSy^{k}\|^{2}-\alpha^{k}\beta^{k}g_{1}(\|JTz^{k}-Jx^{k}\|)\\
&\leq\alpha^{k}\phi(u,x^{k})+\beta^{k}\phi(u,z^{k})+\gamma^{k}\phi(u,y^{k})-\alpha^{k}\beta^{k}g_{1}(\|JTz^{k}-Jx^{k}\|)\\
&\leq\phi(u,x^{k})-\alpha^{k}\beta^{k}g_{1}(\|JTz^{k}-Jx^{k}\|),
\end{aligned}
\end{equation*}
and by similar argument, there exists a continuous, strictly
increasing and convex function\linebreak
 $g_{2}:[0,2r_{2}]\rightarrow\mathbb{R}$ with $g_{2}(0)=0$ such that for $u\in \Omega$, we get
\begin{equation*}
\phi(u,x^{k+1})\leq\phi(u,x^{k})-\alpha^{k}\gamma^{k}g_{2}(\|JSy^{k}-Jx^{k}\|),
\end{equation*}
which imply
\begin{equation}\label{eq16}
\alpha^{k}\beta^{k}g_{1}(\|JTz^{k}-Jx^{k}\|)\leq\phi(u,x^{k})-\phi(u,x^{k+1}),
\end{equation}
\begin{equation}\label{eq17}
\alpha^{k}\gamma^{k}g_{2}(\|JSy^{k}-Jx^{k}\|)\leq\phi(u,x^{k})-\phi(u,x^{k+1}).
\end{equation}
Taking the limits as $k\rightarrow \infty$ in inequalities
(\ref{eq16}) and (\ref{eq17}), we have
\begin{equation}
\lim\limits_{k\rightarrow\infty}g_{1}(\|JTz^{k}-Jx^{k}\|)=0\quad
\& \quad
\lim\limits_{k\rightarrow\infty}g_{2}(\|JSy^{k}-Jx^{k}\|)=0.
\end{equation}
From the properties of $g_{1}$ and $g_{2}$, we get
\begin{equation}\label{eq50}
\lim\limits_{k\rightarrow\infty}\|JTz^{k}-Jx^{k}\|=0\quad \& \quad
\lim\limits_{k\rightarrow\infty}\|JSy^{k}-Jx^{k}\|=0.
\end{equation}
Also, from Lemma \ref{eq6}, (\ref{eq50}) and inequality (\ref{eq500}), we have
\begin{equation*}
\begin{aligned}
\phi(x^{k},x^{k+1})&\leq\|x^{k}\|\|Jx^{k}-(\alpha^{k}Jx^{k}+\beta^{k}JTz^{k}+\gamma^{k}JSy^{k})\|\\
&\quad+\|J^{-1}(Jx^{k})-J^{-1}(\alpha^{k}Jx^{k}+\beta^{k}JTz^{k}+\gamma^{k}JSy^{k})\|\|\alpha^{k}Jx^{k}+\beta^{k}JTz^{k}+\gamma^{k}JSy^{k}\|\\
&\leq(\|x^{k}\|+\frac{2}{c^{2}}\|\alpha^{k}Jx^{k}+\beta^{k}JTz^{k}+\gamma^{k}JSy^{k}\|)\|(Jx^{k})-(\alpha^{k}Jx^{k}+\beta^{k}JTz^{k}+\gamma^{k}JSy^{k})\|\\
&\leq(\|x^{k}\|+\frac{2}{c^{2}}\|\alpha^{k}Jx^{k}+\beta^{k}JTz^{k}+\gamma^{k}JSy^{k}\|)\\
&\quad\times(\alpha^{k}\|Jx^{k}-Jx^{k}\|+\beta^{k}\|Jx^{k}-JTz^{k}\|+\gamma^{k}\|Jx^{k}-JSy^{k}\|)\\
&\rightarrow0 \ \ \text{as}\ \ k\rightarrow\infty,
\end{aligned}
\end{equation*}
consequently, by Lemma \ref{eq22}, we obtain $\lim\limits_{k\rightarrow\infty}\|x^{k}-x^{k+1}\|=0.$ So, $\{x^{k}\}$ converges strongly to $p\in C$. Since $J^{-1}$ is
uniformly norm-to-norm continuous on bounded sets, from
(\ref{eq18}) and (\ref{eq50}), we obtain
\begin{equation}\label{eq29}
\lim\limits_{k\rightarrow\infty}\|Tz^{k}-x^{k}\|=\lim\limits_{k\rightarrow\infty}\|J^{-1}(JTz^{k})-J^{-1}(Jx^{k})\|=0,
\end{equation}
\begin{equation}\label{eq45}
\lim\limits_{k\rightarrow\infty}\|Sy^{k}-x^{k}\|=\lim\limits_{k\rightarrow\infty}\|J^{-1}(JSy^{k})-J^{-1}(Jx^{k})\|=0.
\end{equation}
Combining inequalities (\ref{eq10}) and (\ref{eq20}), we get
\begin{equation*}
\begin{aligned}
\phi(u,x^{k+1})&\leq\alpha^{k}\phi(u,x^{k})+\beta^{k}\phi(u,z^{k})+\gamma^{k}\phi(u,y^{k})\\
&\leq\alpha^{k}\phi(u,x^{k})+\beta^{k}\phi(u,x^{k})+\gamma^{k}\phi(u,y^{k})\\
&=(1-\gamma^{k})\phi(u,x^{k})+\gamma^{k}\phi(u,y^{k}),\\
&\leq(1-\gamma^{k})\phi(u,x^{k})+\gamma^{k}(\phi(u,x^{k})+2\tau(\frac{2\tau}{c^{2}}-\alpha)(\|Ax^{k}-Au\|^{2})\\
&=\phi(u,x^{k})+2\tau\gamma^{k}(\frac{2\tau}{c^{2}}-\alpha)(\|Ax^{k}-Au\|^{2}),
\end{aligned}
\end{equation*}
also, combining inequalities (\ref{eq21}) and (\ref{eq20}), we
have
\begin{equation*}
\begin{aligned}
\phi(u,x^{k+1})&\leq\alpha^{k}\phi(u,x^{k})+\gamma^{k}\phi(u,y^{k})+\beta^{k}\phi(u,z^{k})\\
&\leq(1-\beta^{k})\phi(u,x^{k})+\beta^{k}\phi(u,z^{k})\\
&=\phi(u,x^{k})+2\tau\beta^{k}(\frac{2\tau}{c^{2}}-\alpha)(\|Au^{k}-Au\|^{2}).
\end{aligned}
\end{equation*}
Therefore, we get
\begin{equation*}
2\tau\gamma^{k}(\alpha-\frac{2\tau}{c^{2}})(\|Ax^{k}-Au\|^{2})\leq\phi(u,x^{k})-\phi(u,x^{k+1}),
\end{equation*}
\begin{equation*}
2\tau\beta^{k}(\alpha-\frac{2\tau}{c^{2}})(\|Au^{k}-Au\|^{2})\leq\phi(u,x^{k})-\phi(u,x^{k+1}).
\end{equation*}
Since $\{\phi(u,x^{k})\}$ is convergent, it follows from
conditions (ii) and (iv) that
\begin{equation}\label{eq.m.1}
\lim\limits_{k\rightarrow\infty}\|Ax^{k}-Au\|^{2}=0\quad \& \quad
\lim\limits_{k\rightarrow\infty}\|Au^{k}-Au\|^{2}=0.
\end{equation}
From (\ref{eq5}), (\ref{eq.m.1}), lemmas \ref{eq4} and \ref{eq6}
and assumption
 $\Vert Ax\Vert\leq\Vert Ax-Au\Vert$ for all $x\in C$ and all $u\in  \Omega$, we obtain
\begin{equation*}
\begin{aligned}
\phi(x^{k},y^{k})&=\phi(x^{k},\Pi_{C}J^{-1}(Jx^{k}-\tau Ax^{k}))\\
&\leq\phi(x^{k},J^{-1}(Jx^{k}-\tau Ax^{k}))\\
&=V(x^{k},Jx^{k}-\tau Ax^{k})\\
&\leq V(x^{k},Jx^{k}-\tau Ax^{k}+\tau Ax^{k})-2\langle J^{-1}(Jx^{k}-\tau Ax^{k})-x^{k},\tau Ax^{k}\rangle\\
&=\phi(x^{k},x^{k})+2\langle J^{-1}(Jx^{k}-\tau Ax^{k})-x^{k},-\tau Ax^{k}\rangle\\
&=2\langle J^{-1}(Jx^{k}-\tau Ax^{k})-J^{-1}J(x^{k}),-\tau Ax^{k}\rangle\\
&\leq\|J^{-1}(Jx^{k}-\tau Ax^{k})-J^{-1}(Jx^{k})\|\|\tau Ax^{k}\|\\
&\leq\frac{4}{c^{2}}\|JJ^{-1}(Jx^{k}-\tau Ax^{k})-JJ^{-1}(Jx^{k})\|\|\tau Ax^{k}\|\\
&=\frac{4}{c^{2}}\tau^{2}\|Ax^{k}\|^{2}\\
&\leq\frac{4}{c^{2}}\tau^{2}\|Ax^{k}-Au\|^{2}\\
&\rightarrow0 \ \ \text{as}\ \ k\rightarrow\infty,
\end{aligned}
\end{equation*}
and in the same way, we can conclude that
\begin{equation*}
\phi(u^{k},z^{k})=\phi(u^{k},\Pi_{C}J^{-1}(Ju^{k}-\tau Au^{k}))\leq
\frac{4}{c^{2}}\tau^{2}\|Au^{k}-Au\|^{2}\rightarrow0 \ \
\text{as}\ \ k\rightarrow\infty.
\end{equation*}
Consequently by Lemma \ref{eq22}, we obtain
\begin{equation}\label{eq27}
\lim\limits_{k\rightarrow\infty}\|x^{k}-y^{k}\|=\lim\limits_{k\rightarrow\infty}\|u^{k}-z^{k}\|=0.
\end{equation}
\par
Let $r_{3}=\sup_{k\geq1}\{\|u^{k}\|,\|x^{k}\|\}$. So, by Lemma
\ref{kamimra}, there exists a continuous, strictly increasing and
convex function $g_{3}:[0,2r_{3}]\rightarrow\mathbb{R}$ with
$g_{3}(0)=0$ such that
$$g_{3}(\|u^{k}-x^{k}\|)\leq\phi(u^{k},x^{k}).$$
Let $u\in \Omega$, since $\phi(u,Tz^{k})\leq\phi(u,u^{k})$ and
$u^{k}=K_{r^{k}}x^{k}$, we observe from condition (v) of Lemma
\ref{eq13} that
\begin{equation*}
\begin{aligned}
g_{3}(\|u^{k}-x^{k}\|)&\leq\phi(u^{k},x^{k})\\
&\leq\phi(u,x^{k})-\phi(u,u^{k})\\
&\leq\phi(u,x^{k})-\phi(u,Tz^{k})\\
&=\|u\|^{2}+\|x^{k}\|^{2}-2\langle u,Jx^{k}\rangle-\|u\|^{2}-\|Tz^{k}\|^{2}+2\langle u,JTz^{k}\rangle\\
&=\|x^{k}\|^{2}-\|Tz^{k}\|^{2}+2\langle u,JTz^{k}-Jx^{k}\rangle\\
&\leq\|x^{k}\|^{2}-\|Tz^{k}\|^{2}+2\|u\|\|JTz^{k}-Jx^{k}\|\\
&\leq(\|x^{k}-Tz^{k}\|+\|Tz^{k}\|)^{2}-\|Tz^{k}\|^{2}+2\|u\|\|JTz^{k}-Jx^{k}\|\\
&=\|x^{k}-Tz^{k}\|^{2}+2\|x^{k}-Tz^{k}\|\|Tz^{k}\|+2\|u\|\|JTz^{k}-Jx^{k}\|.
\end{aligned}
\end{equation*}
From equalities (\ref{eq50}) and (\ref{eq29}), we have
$\lim\limits_{k\rightarrow\infty}g_{3}(\|u^{k}-x^{k}\|)=0$ and so
\begin{equation}\label{eq48}
\lim\limits_{k\rightarrow\infty}\|u^{k}-x^{k}\|=0.
\end{equation}
Therefore, from (\ref{eq27}) and (\ref{eq48}), we have
\begin{equation}\label{eq28}
\|x^{k}-z^{k}\|\leq\|x^{k}-u^{k}\|+\|u^{k}-z^{k}\| \rightarrow0 \
\ \text{as}\ \ k\rightarrow\infty.
\end{equation}
It follows from (\ref{eq29}), (\ref{eq45}), (\ref{eq27}) and
(\ref{eq28}) that
\begin{equation}
\|Tz^{k}-z^{k}\|\leq\|Tz^{k}-x^{k}\|+\|x^{k}-z^{k}\|\rightarrow0,\
\ \text{as} \ \ k\rightarrow\infty,
\end{equation}
\begin{equation}
\|Sy^{k}-y^{k}\|\leq\|Sy^{k}-x^{k}\|+\|x^{k}-y^{k}\|\rightarrow0,\
\ \text{as} \ \ k\rightarrow\infty.
\end{equation}
 From (\ref{eq27}) and (\ref{eq28}), we can conclude that $\{y^{k}\}$ and $\{z^{k}\}$ converge strongly to $p\in C$, using the definitions of $T$ and $\hat{F}(T)$, we have $p\in \hat{F}(T)=F(T)$. Also the definitions of $S$ and $\hat{F}(S)$ imply that $p\in \hat{F}(S)=F(S)$. Hence, $p\in F(T)\cap F(S)$.
\par
Now, we show that $p\in GEP(f,A)$. From (\ref{eq27}) and
(\ref{eq30}), we obtain
\begin{equation}
\lim\limits_{k\rightarrow\infty}\|Ju^{k}-Jx^{k}\|=0,
\end{equation}
since $J$ is uniformly norm-to-norm continuous on bounded sets. It
follows from condition (iii) that
$\lim\limits_{k\rightarrow\infty}\frac{\|Ju^{k}-Jx^{k}\|}{r^{k}}=0.$
By the definition of $u^{k}=K_{r^{k}}x^{k}$, we get
$F(u^{k},y)+\frac{1}{r^{k}}\langle
y-u^{k},Ju^{k}-Jx^{k}\rangle\geq0$, for all $y\in C$, where
$F(u^{k},y)=f(u^{k},y)+\langle Au^{k},y-u^{k}\rangle$. It is
easily seen that $y\rightarrow f(x,y)+\langle Ax,y-x\rangle$ is
convex and lower semicontinuous, so from Lemma \ref{eq23}, it is
weakly lower semicontinuous. Thus bifunction $F:C\times
C\rightarrow\mathbb{R}$ satisfying the condition $(A4)$ and
clearly satisfying in $(A1)-(A3)$. We have from $(A2)$ that
\begin{equation*}
\frac{1}{r^{k}}\langle y-u^{k},Ju^{k}-Jx^{k}\rangle\geq
-F(u^{k},y)\geq F(y,u^{k}),
\end{equation*}
for all $y\in C$. Taking the limit inferior on both sides of the
last inequality and using $(A4)$, we can conclude that
$$F(y,p)\leq0, \qquad \forall\: y\in C.$$ Let $y_{t}=ty+(1-t)p$ for
all $y\in C$ and all $0<t<1$, the convexity of $C$ implies that
$y_{t}\in C$ and hence $F(y_{t},p)\leq0$. Therefore, from $(A1)$
and $(A4)$ we have
$$0=F(y_{t},y_{t})\leq tF(y_{t},y)+(1-t)F(y_{t},p)\leq tF(y_{t},y).$$
Hence, $F(y_{t},y)\geq0$ for all $y\in C$. Taking the limit as
$t\downarrow0$ and using $(A3)$, we yield that $F(p,y)\geq0$ and
therefore $f(p,y)+\langle Ap,y-p\rangle\geq0$ for all $y\in C$, so
$p\in GEP(f,A)$.
\par
Now, we prove that $p\in SOL(C,A)$. Let $B\subset E\times E^{*}$
be an operator which is defined as follows:
$$B\nu=\begin{cases}A\nu+N_{C}(\nu), \ \ \nu\in C,\\
\emptyset, \ \ \ \ \ \ \ \ \ \ \ \ \  \ \ \nu\not \in C,
\end{cases}$$
also, we know that $A$ is hemicontinuous, because
\begin{equation}\label{eq501}
\lambda\langle x-y,A(\lambda x+(1-\lambda)y)-Ay\rangle\geq\alpha\|A(\lambda x+(1-\lambda)y)-Ay\|^{2},
\end{equation}
for all $0\leq\lambda\leq1$. Taking the limits as $\lambda\rightarrow 0$ in inequality (\ref{eq501}), therefore\linebreak $A(\lambda x+(1-\lambda)y)\rightarrow Ay.$ So, it follows from Lemma \ref{eq34} that $B$ is maximal monotone and
$B^{-1}(0)= SOL(C,A)$.
 Let $(\nu,w)\in G(B)$. Since $w\in B\nu=A\nu+N_{C}(\nu)$, we get $w-A\nu\in N_{C}(\nu)$. Since $y^{k}\in C$, we obtain
\begin{equation}\label{eq46}
\langle \nu-y^{k},w-A\nu\rangle\geq0.
\end{equation}
From Lemma \ref{eq4} and (\ref{eq36}), we get
\begin{equation*}
\begin{aligned}
2\langle\nu-y^{k},Jy^{k}-J(J^{-1}(Jx^{k}-\tau Ax^{k}))\rangle &=\phi(\nu,J^{-1}(Jx^{k}-\tau Ax^{k}))-\phi(\nu,y^{k})\\
&\quad-\phi(y^{k},J^{-1}(Jx^{k}-\tau Ax^{k}))\\
&\geq 0.
\end{aligned}
\end{equation*}
Thus, $\: \langle\nu-y^{k},Jy^{k}-J(J^{-1}(Jx^{k}-\tau
Ax^{k})\rangle\geq0$. Hence
\begin{equation}\label{eq57}
\: \langle\nu-y^{k},Ax^{k}+\frac{Jy^{k}-Jx^{k}}{\tau}\rangle\geq0.
\end{equation}
Using the definition of $A$ and inequalities (\ref{eq46}) and
(\ref{eq57}), we have
\begin{equation*}
\begin{aligned}
\langle\nu-y^{k},w\rangle &\geq \langle\nu-y^{k},A\nu \rangle\\
&\geq\langle\nu-y^{k},A\nu\rangle-\langle\nu-y^{k},\frac{Jy^{k}-Jx^{k}}{\tau}+Ax^{k}\rangle\\
&=\langle\nu-y^{k},A\nu-Ay^{k} \rangle+\langle\nu-y^{k},Ay^{k}-Ax^{k} \rangle\\
&\quad-\langle\nu-y^{k},\frac{Jy^{k}-Jx^{k}}{\tau}\rangle\\
&\geq-\langle\nu-y^{k},Ax^{k}-Ay^{k}\rangle-\langle\nu-y^{k},\frac{Jy^{k}-Jx^{k}}{\tau}\rangle\\
&\geq-(\| Ax^{k}-Ay^{k}\|+\frac{1}{\tau}\| Jy^{k}-Jx^{k}\|)\|\nu-y^{k}\|\\
&\geq-(\frac{1}{\alpha}\|x^{k}-y^{k}\|+\frac{1}{\tau}\|Jy^{k}-Jx^{k}\|)\|\nu-y^{k}\|.
\end{aligned}
\end{equation*}
By letting $k\rightarrow\infty$ and using (\ref{eq30}) and
(\ref{eq27}), we obtain $\langle \nu-p,w\rangle\geq0,$ therefore
$p\in B^{-1}(0)=SOL(C,A)$, because of $B$ is a maximal monotone
operator.
\par
Now, let $\nu^{k}=\Pi_{\Omega}(x^{k})$, from  inequality
(\ref{eq20}), we have
\begin{equation}\label{eq39}
\phi(\nu^{k},x^{k+1})\leq\phi(\nu^{k},x^{k}),
\end{equation}
hence, from Lemma \ref{eq39}, we get
\begin{equation*}
\phi(\nu^{k+1},x^{k+1})=\phi(\Pi_{\Omega}(x^{k+1}),x^{k+1})\leq\phi(\nu^{k},x^{k+1})\leq\phi(\nu^{k},x^{k}).
\end{equation*}
This implies that
$\lim\limits_{k\rightarrow\infty}\phi(\nu^{k},x^{k})$ exists. So,
$\{\phi(\nu^{k},x^{k})\}$ is bounded. Using inequality
(\ref{eq14}), we conclude that $\{\nu^{k}\}$ is bounded. Since
$\nu^{k+m}=\Pi_{\Omega}(x^{k+m}),$ for all $m\in \mathbb{N}$, from
Lemma \ref{eq4} and inequality (\ref{eq39}), we obtain
\begin{equation*}
\phi(\nu^{k},\nu^{k+m})+\phi(\nu^{k+m},x^{k+m})\leq\phi(\nu^{k},x^{k+m})\leq\phi(\nu^{k},x^{k}).
\end{equation*}
So,
\begin{equation*}
\phi(\nu^{k},\nu^{k+m})\leq\phi(\nu^{k},x^{k})-\phi(\nu^{k+m},x^{k+m}).
\end{equation*}
\par
Let $\acute{r}=\sup_{k\geq1}\|\nu^{k}\|$. Using the Lemma
\ref{kamimra}, there exists a continuous strictly increasing and
convex function $\acute{g}$ with $\acute{g}(0)=0$ such that
\begin{equation*}
\begin{aligned}
\acute{g}(\|\nu^{k}-\nu^{k+m}\|)&\leq\phi(\nu^{k},\nu^{k+m})\leq\phi(\nu^{k},x^{k})-\phi(\nu^{k+m},x^{k+m}).
\end{aligned}
\end{equation*}
Since $\lim\limits_{k\rightarrow\infty}\phi(\nu^{k},x^{k})$
exists, from the properties $\acute{g}$, we have $\{\nu^{k}\}\in
\Omega$ is a cauchy sequence. Since $\Omega$ is closed, so
$\{\nu^{k}\}$ converges strongly to $u^{*}\in \Omega$ and from
Lemma \ref{eq4}, we get
$\langle\nu^{k}-x^{k},Jp-J\nu^{k}\rangle\geq0$. Therefore, we get
$\langle u^{*}-p,Jp-Ju^{*}\rangle\geq0$. On the other hand, since
$J$ is monotone, so $\langle u^{*}-p,Jp-Ju^{*}\rangle\leq0$. Thus
$\langle u^{*}-p,Jp-Ju^{*}\rangle=0$, since $J$ is one-one, we get
$p=u^{*}$. Therefore $x^{k}\rightarrow u^{*}$ and inequalities
(\ref{eq27}) and (\ref{eq28}) imply that $y^{k}\rightarrow u^{*}$
and $z^{k}\rightarrow u^{*}$, where
$u^{*}=\lim\limits_{k\rightarrow\infty}\Pi_{\Omega}(x^{k})$. \qed
\begin{cor}
Let $C$ be a nonempty closed convex subset of a $2$-uniformly
convex, uniformly smooth Banach space $E$. Let $A:C\rightarrow
E^{*}$ ba an $\alpha $-inverse strongly monotone operator and
$S:C\rightarrow C$ be a relatively nonexpansive mapping  such that
$\Omega:=SOL(C,A)\cap F(S)\neq\emptyset$ and $\Vert
Ax\Vert\leq\Vert Ax-Au\Vert   \  \text{for all}\ x\in C\ \text{and
all} \ u\in  \Omega$. Assume that $\Pi_{C}$ is the generalized
metric projection from $E$ onto $C$. Let $\{x^{k}\}$ be a
sequence generated by $x^{1}\in C$ and
\begin{equation}\label{eq24}
\begin{cases}
 u^{k}\in C\;\text{such that}\quad \langle Au^{k},y-u^{k}\rangle +\frac{1}{r^{k}}\langle y-u^{k},Ju^{k}-Jx^{k}\rangle\geq0,\quad\forall\: y\in C,&\\
  y^{k}:=\Pi_{C}J^{-1}(Jx^{k}-\tau Ax^{k}),&\\
 z^{k}:=\Pi_{C}J^{-1}(Ju^{k}-\tau Au^{k}),&\\
 x^{k+1}:=\Pi_{C}J^{-1}(\alpha^{k}Jx^{k}+\beta^{k}Jz^{k}+\gamma^{k}JSy^{k}).&
  \end{cases}
\end{equation}
Furthermore, assume that $\{\alpha^{k}\}$, $\{\beta^{k}\}$ and
$\{\gamma^{k}\}$ are three sequences in $[0,1]$ satisfying the
following conditions:
 \begin{enumerate}
      \item[(i)]$\alpha^{k}+\beta^{k}+\gamma^{k}=1$,
      \item[(ii)]$\liminf\limits_{k\rightarrow \infty} \alpha^{k}\beta^{k}>0 \quad \& \quad \liminf\limits_{k\rightarrow \infty} \alpha^{k}\gamma^{k}>0$;
      \item[(iii)]$\{r^{k}\}\subset [a,\infty)$ for some $a>0$;
      \item[(iv)] $0<\tau<\frac{c^{2}\alpha}{2}$, where $\frac{1}{c}\,\,(0< c\leq 1)$ is the $2$-uniformly convexity constant of $E$.
    \end{enumerate}
Then the sequences $\{x^{k}\}_{k=1}^{\infty}$,
$\{y^{k}\}_{k=1}^{\infty}$ and $\{z^{k}\}_{k=1}^{\infty}$
generated by (\ref{eq24}) converge strongly  to the some solution
$u^{*}\in \Omega$, where $u^{*}=\lim\limits_{k\rightarrow
\infty}\Pi_{\Omega}({x^{k}}).$
\end{cor}
{\it Proof} Letting $f\equiv0$ and $T=I$ in Theorem \ref{eq25}, we
get the desired result.\qed
\section{Numerical example}
Now, we demonstrate Theorem \ref{eq25} with an example.
\begin{exa}
Let $E=\mathbb{R}$, $C=[-4,4]$ and $A=I$ such that $\alpha=1$,
$c=1$ and $\tau=\frac{1}{4}.$ Define $f(u,y):=9y^{2}+3uy-12u^{2},$
 we see that $f$ satisfies the conditions $(A1)-(A4)$ as follows:
 \begin{description}
 \item[(A1)] $f(u,u):=9u^{2}+3u^{2}-12u^{2}=0$ for all $u\in [-4,4]$,
 \item[(A2)] $f(u,y)+f(y,u)=-3(y-u)^{2}\leq0$ for all $u, y\in [-4,4]$, i.e., $f$ is monotone,
 \item[(A3)] for each $u,y,z\in [-4,4]$,
 \begin{equation*}
 \begin{aligned}
 \lim\limits_{t\downarrow0}f(tz+(1-t)u,y)&=\lim\limits_{t\downarrow0}(9y^{2}+3(tz+(1-t)u)y-12(tz+(1-t)u)^{2})\\
 &=9y^{2}+3uy-12u^{2}\\
 &=f(u,y),
 \end{aligned}
 \end{equation*}
 \item[(A4)] it is easily seen that for each $u\in [-4,4]$, $y\rightarrow (9y^{2}+3uy-12u^{2})$
  is convex and lower semicontinuous.
\end{description}
\par
On the other hand, we have $\langle Au, y-u\rangle=\langle u,
y-u\rangle=u(y-u)=uy-u^{2}.$ Also
 $$\frac{1}{r}\langle y-u,Ju-Jx\rangle=\frac{1}{r}(y-u)(u-x)=\frac{1}{r}(uy-u^{2}+ux-xy).$$
From condition (i) of Lemma \ref{eq13}, $K_{r}$ is single-valued,
let $u=K_{r}x$, for any $y\in [-4,4]$ and $r>0$, we have
$$f(u,y)+\langle Au,y-u\rangle+\frac{1}{r}\langle y-u,Ju-Jx\rangle\geq0.$$
Thus,
\begin{equation*}
\begin{aligned}
9ry^{2}+3ruy&-12ru^{2}+ruy-ru^{2}+uy-u^{2}+ux-xy\\
&=9ry^{2}+(4ru+u-x)y-13ru^{2}-u^{2}+ux\\
&\geq0.
\end{aligned}
\end{equation*}
Now, let $a=9r$, $b=4ru+u-x$ and $c=-13ru^{2}-u^{2}+ux$.\\
Hence, we should have $\Delta=b^{2}-4ac\leq0$, i.e.,
\begin{equation*}
\begin{aligned}
\Delta&=((4r+1)u-x)^{2}+36ru((13r+1)u-x)\\
&=484r^{2}u^{2}+44ru^{2}+u^{2}+x^{2}-44rux-2ux\\
&=((22r+1)u-x)^{2}\\
&\leq0.
\end{aligned}
\end{equation*}
So, it follows that $u=\frac{x}{1+22r}$. Therefore,
$K_{r}x=\frac{x}{1+22r}$. So, using the notations of Theorem
\ref{eq25}, we have
$u^{k}=K_{r^{k}}x^{k}=\frac{x^{k}}{1+22r^{k}}$.
Since $F(K_{r^{k}})=0$, from condition (iii) of Lemma \ref{eq13}, $GEP(f,I)=0$.\\
Define $T:C\rightarrow C$ by $Tx=x$ for all $x\in C$, therefore
$F(T)=C$ and $\phi(p,Tx)=\phi(p,x)$, for all $x\in C$ and all
$p\in F(T)$. Let $x^{k}\rightharpoonup p,$ so
$\lim\limits_{k\rightarrow\infty}(Tx^{k}-x^{k})=0$, these imply
that $\hat{F}(T)=C$. Therefore, $\hat{F}(T)=F(T)$, i.e., $T$ is
relatively nonexpansive mapping.
\par
 Now, define $S:C\rightarrow C$ by
$Sx=\frac{2}{9}x$ for all $x\in C$, so $F(S)=\{0\}$ and
\begin{equation*}
\begin{aligned}
\phi(0,Sx)&=\phi(0,\frac{2}{9}x)\\
&=0-2\langle 0,\frac{2}{9}x\rangle+|\frac{2}{9}x|^{2}\\
&\leq|x|^{2}\\
&=\phi(0,x),
\end{aligned}
\end{equation*}
for all $x\in C$.
 \begin{figure}[!h]
\includegraphics[totalheight=2.9in]{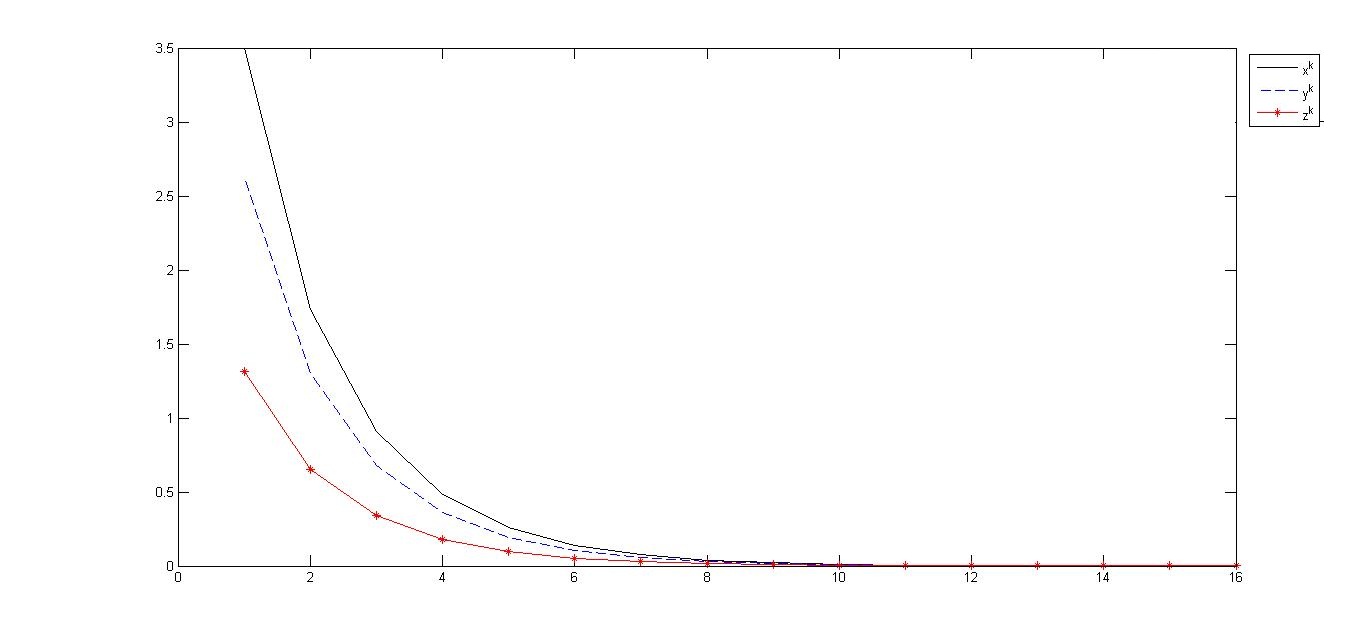}\vspace{-.55cm}\caption{}
\end{figure}
\begin{figure}[!h]
\includegraphics[totalheight=2.9in]{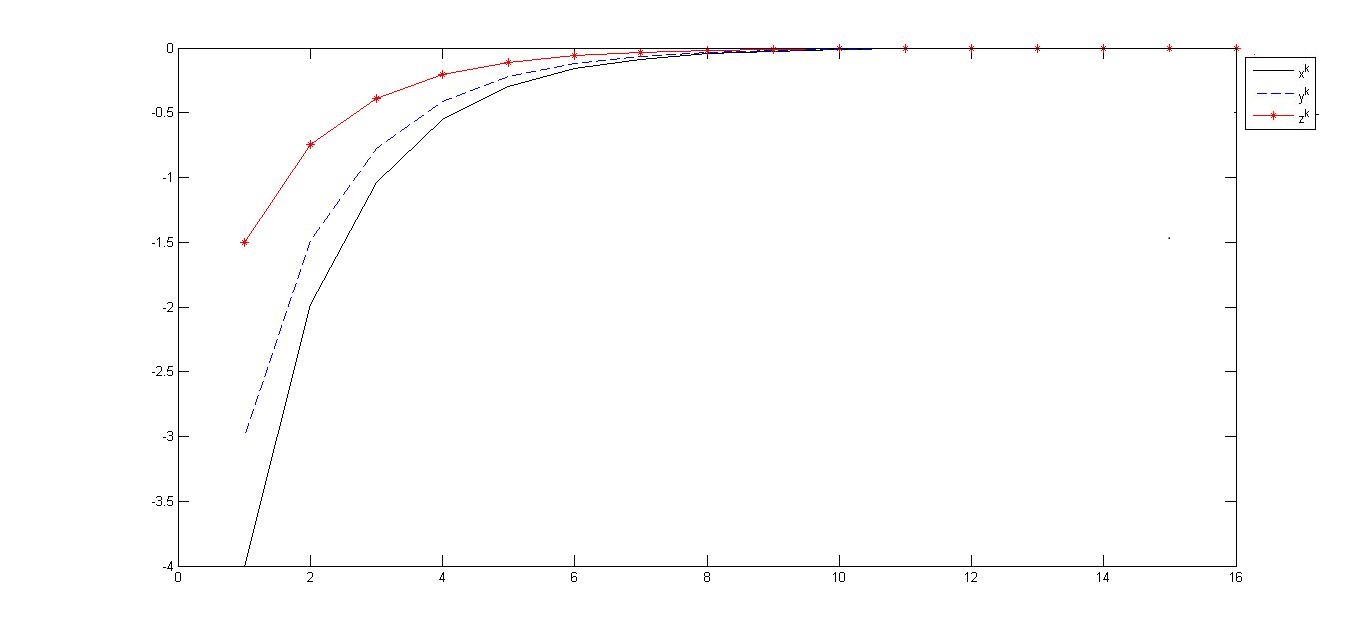}\vspace{-.55cm}\caption{}
\end{figure}
It is easy to see that $S$ is also relatively nonexpansive
mapping. Moreover, since $SOL(C,I)=\{u\in C; \langle
u,y-u\rangle\geq0\}=\{0\}$. So
$\Omega=\{0\}$.
\par
 Assume that $\alpha^{k}=\frac{1}{3}+\frac{1}{4k}$,
$\beta^{k}=\frac{1}{2}-\frac{1}{6k}$,
$\gamma^{k}=\frac{1}{6}-\frac{1}{12k}$ and $r^{k}=\frac{1}{22}$,
for all $k \geq 1,$ so $\{\alpha^{k}\}$, $\{\beta^{k}\}$ and
$\{\gamma^{k}\}$ are three sequences in $[0,1]$ such that
satisfies the conditions (i)-(iii) in the Theorem \ref{eq25}.
 Since $C$ is symmetric and $u^{k}=\frac{1}{2}x^{k}$, we have
\begin{equation}\label{eq72}
\begin{cases}
 y^{k}:=\Pi_{C}(x^{k}-\frac{1}{4}x^{k})=\frac{3}{4}x^{k},&\\
 z^{k}:=\Pi_{C}(u^{k}-\frac{1}{4}u^{k})=\frac{3}{4}u^{k}=\frac{3}{8}x^{k},&
  \end{cases}
\end{equation}
also
\begin{equation}\label{eq71}
\begin{aligned}
 x^{k+1}:&=\Pi_{C}(\alpha^{k}x^{k}+\beta^{k}Tz^{k}+\gamma^{k}Sy^{k})\\
 &=\alpha^{k}x^{k}+\beta^{k}(\frac{3}{8}x^{k})+\gamma^{k}(\frac{1}{6}x^{k})\\
 &=((\frac{1}{3}+\frac{1}{4k})+\frac{3}{8}(\frac{1}{2}-\frac{1}{6k})+\frac{1}{6}(\frac{1}{6}-\frac{1}{12k}))x^{k}\\
 &=(\frac{79}{144}-\frac{16}{304k})x^{k}.
 \end{aligned}
 \end{equation}
 \begin{center}
\vspace{3mm}
\begin{tabular}{ c c c c}
\hline  & \hspace{2 cm} {\small Numerical Results for $x^{1}
=3.5$}&
\\
\cline{1-4}
k\hspace{-.3 cm}  & \hspace{-1.9 cm} {$x^{k}$} &\hspace{-3.8 cm}{$y^{k}$} &\hspace{0 cm}{$z^{k}$}\\
\hline\hline
1\hspace{-.3 cm} &\hspace{-1.7 cm}$3.5$&\hspace{-3.8 cm} $2.625$ & \hspace{0 cm} $1.3125$ \\
2\hspace{-.3 cm} &\hspace{-1.7 cm}$1.7359 $&\hspace{-3.8 cm}$1.3019$ & \hspace{0 cm}$0.65$\\
3\hspace{-.3 cm} &\hspace{-1.7 cm}$0.9067 $&\hspace{-3.8 cm}$0.68$ & \hspace{0 cm}$0.34$\\
\hspace{-.3 cm} &\hspace{-1.7 cm}$\vdots$ & \hspace{-3.8 cm}$\vdots$ & \hspace{0 cm}$\vdots$\\
45\hspace{-.3 cm} &\hspace{-1.7 cm}$7.69e-12$&\hspace{-3.8 cm}$5.77e-12$ & \hspace{0 cm}$2.89e-12$\\
46\hspace{-.3 cm} &\hspace{-1.7 cm}$4.21e-12 $&\hspace{-3.8 cm}$3.16e-12$ & \hspace{0 cm}$1.58e-12$\\
47\hspace{-.3 cm} &\hspace{-1.7 cm}$2.31e-12$&\hspace{-3.8 cm}$1.73e-12$ & \hspace{0 cm}$8.65e-13$\\
\hspace{-.3 cm} &\hspace{-1.7 cm}$\vdots $&\hspace{-3.8cm}$\vdots$ & \hspace{0 cm}$\vdots$\\
98\hspace{-.3 cm} &\hspace{-1.7 cm}$5.93e-26$&\hspace{-3.8 cm}$4.45e-26$ & \hspace{0 cm}$2.22e-26$\\
99\hspace{-.3 cm} &\hspace{-1.7 cm}$3.25e-26$&\hspace{-3.8 cm}$2.44e-26$ & \hspace{0 cm}$1.22e-26$\\
100\hspace{-.3 cm} &\hspace{-1.7 cm}$1.78e-26$&\hspace{-3.8 cm}$1.34e-26$ & \hspace{0 cm}$6.68e-27$\\
\hline
\end{tabular}
\end{center}
\vspace{1 cm}
 \begin{center}
\vspace{3mm}
\begin{tabular}{ c c c c}
\hline  & \hspace{2 cm} {\small Numerical Results for $x^{1}
=-4$}&
\\
\cline{1-4}
k\hspace{-.3 cm}  & \hspace{-1.9 cm} {$x^{k}$} &\hspace{-3.8 cm}{$y^{k}$} &\hspace{0 cm}{$z^{k}$}\\
\hline\hline
1\hspace{-.3 cm} &\hspace{-1.7 cm}$-4 $&\hspace{-3.8 cm} $-3$ & \hspace{0 cm} $-1.5$ \\
2\hspace{-.3 cm} &\hspace{-1.7 cm}$-1.9392$&\hspace{-3.8 cm}$-1.488$ & \hspace{0 cm}$-0.744$\\
3\hspace{-.3 cm} &\hspace{-1.7 cm}$-1.03619$&\hspace{-3.8 cm}$-0.777$ & \hspace{0 cm}$-0.389$\\
\hspace{-.3 cm} &\hspace{-1.7 cm}$\vdots$ & \hspace{-3.8 cm}$\vdots$ & \hspace{0 cm}$\vdots$\\
45\hspace{-.3 cm} &\hspace{-1.7 cm}$-8.79e-12 $&\hspace{-3.8 cm}$-6.59e-12$ & \hspace{0 cm}$-3.3e-12$\\
46\hspace{-.3 cm} &\hspace{-1.7 cm}$-4.81e-12 $&\hspace{-3.8 cm}$-3.61e-12$ & \hspace{0 cm}$-1.8e-12$\\
47\hspace{-.3 cm} &\hspace{-1.7 cm}$-2.64e-12 $&\hspace{-3.8 cm}$-1.98e-12$ & \hspace{0 cm}$-9.88e-13$\\
\hspace{-.3 cm} &\hspace{-1.7 cm}$\vdots $&\hspace{-3.8 cm}$\vdots$ & \hspace{0 cm}$\vdots$\\
98\hspace{-.3 cm} &\hspace{-1.7 cm}$-6.78e-26 $&\hspace{-3.8 cm}$-5.09e-26$ & \hspace{0 cm}$-2.54e-26$\\
99\hspace{-.3 cm} &\hspace{-1.7 cm}$-3.72e-26 $&\hspace{-3.8 cm}$-2.79e-26$ & \hspace{0 cm}$-1.39e-26$\\
100\hspace{-.3 cm} &\hspace{-1.7 cm}$-2.04e-26 $&\hspace{-3.8 cm}$-1.53e-26$ & \hspace{0 cm}$-7.64e-27$\\
\hline
\end{tabular}
\end{center}
\vspace{1 cm}
 Since $\Omega=\{0\}$,
 we get $\Pi_{\Omega}(x^{k})=0$ for all $k\geq1$.
 Taking the limit as $k\rightarrow\infty$ in (\ref{eq71}),
 we obtain $\lim\limits_{k\rightarrow\infty}x^{k}=0$ and therefore (\ref{eq72})
 implies that $\lim\limits_{k\rightarrow\infty}y^{k}=\lim\limits_{k\rightarrow\infty}z^{k}=0$.
 See Figure1 and Figure2 for the values $x^{1}=3.5$ and $x^{1}=-4$. The computations associated
 with example were performed using MATLAB software.
\end{exa}

{\footnotesize}

\end{document}